\documentclass[12pt]{amsart}

\usepackage{comment}

\usepackage{extarrows}

\usepackage{amsmath} 
\usepackage{amsthm} 
\usepackage{amsfonts} 
\usepackage{amssymb} 
\usepackage{stmaryrd} 
\usepackage{mathtools} 

\usepackage[all]{xy} 
\usepackage{tikz-cd} 
\usepackage{array} 
\usepackage{tabu} 
\usepackage{enumerate} 

\usepackage[shortlabels]{enumitem} 

\usepackage{fancyhdr} 

\usepackage[marginratio=1:2,paper=letterpaper,tmargin=52pt,lmargin=74pt,rmargin=74pt,includehead]{geometry} 

\usepackage{imakeidx} 

\usepackage[pagebackref=true
,
]{hyperref} 
\hypersetup{colorlinks=true,linkcolor=blue,filecolor=magenta,urlcolor=cyan,} 
\usepackage{bm} 

\makeindex 

\newcommand{\hyp}[1]{\hyperref[#1]{\ref{#1}}} 

\newcommand{\C}{\mathbb{C}} 
\newcommand{\N}{\mathbb{N}} 
\newcommand{\Q}{\mathbb{Q}} 
\newcommand{\R}{\mathbb{R}} 
\newcommand{\Z}{\mathbb{Z}} 

\newcommand{\calA}{\mathcal{A}}

\newcommand{\calE}{\mathcal{E}}

\newcommand{\calH}{\mathcal{H}}

\newcommand{\calL}{\mathcal{L}}

\newcommand{\logdr}[2]{\Omega^\bullet_{#1}(\log #2)} 

\newcommand{\Gr}{\mathrm{Gr}} 
\newcommand{\id}{\mathrm{id}} 

\newtheorem{thm}{Theorem}[section]
\newtheorem{lem}[thm]{Lemma}
\newtheorem{prop}[thm]{Proposition}
\newtheorem{cor}[thm]{Corollary}

\theoremstyle{definition}
\newtheorem{dfn}[thm]{Definition}
\newtheorem{ex}[thm]{Example}

\newtheorem{remk}[thm]{Remark}

\DeclareMathOperator{\Aut}{Aut}

\DeclareMathOperator{\Tors}{Tors}
\DeclareMathOperator{\Ext}{Ext}
\DeclareMathOperator{\Hom}{Hom}

\DeclareMathOperator{\gen}{gen}

\newcommand{\f}[5]{
\begin{array}{rcl}
#1\colon #2 & \longrightarrow & #3 \\
#4 & \longmapsto & #5 \\
\end{array}
}

\newcommand{\cL}{{\mathcal L}}

\newcommand{\cA}{{\mathcal A}}

\newcommand{\ul}[1]{\underline{#1}}

\newcommand{\ov}[1]{\overline{#1}}

\title[]{Hodge theory on Alexander invariants -- a survey}

\author{Eva Elduque}
\address{Department of Mathematics, University of Michigan-Ann Arbor, 530 Church St, Ann Arbor, MI 48109, USA.}
\email {elduque@umich.edu}\urladdr{http://www-personal.umich.edu/~elduque}

\author{Christian Geske}
\address{Department of Mathematics, Northwestern University, 2033 Sheridan Rd, Evanston, IL 60208, USA.}
\email {christian.geske@northwestern.edu}

\author{Mois\'es Herrad\'on Cueto}
\address{Department of Mathematics, Louisiana State University, 303 Lockett Hall, Baton Rouge, LA 70803, USA.}
\email {moises@lsu.edu}\urladdr{http://www.math.lsu.edu/~moises}

\author{Lauren\c{t}iu Maxim}
\address{Department of Mathematics, University of Wisconsin-Madison, 480 Lincoln Drive, Madison WI 53706-1388, USA.}
\email {maxim@math.wisc.edu}\urladdr{https://www.math.wisc.edu/~maxim/}

\author{Botong Wang}
\address{Department of Mathematics, University of Wisconsin-Madison, 480 Lincoln Drive, Madison WI 53706-1388, USA.}
\email {wang@math.wisc.edu}\urladdr{http://www.math.wisc.edu/~wang/}


\dedicatory{Dedicated to the memory of Prof. \c Stefan Papadima}


\begin{document}

\date{\today}

\maketitle

\begin{abstract} 
We survey recent developments in the study of Hodge theoretic aspects of Alexander-type invariants associated with smooth complex algebraic varieties.\end{abstract}


\setlength{\headheight}{15pt}

\thispagestyle{fancy}
\renewcommand{\headrulewidth}{0pt}

\section{Introduction}
In this expository note, we survey recent developments in the study of Hodge theoretic aspects of Alexander-type invariants of complex algebraic manifolds. Our main goal is to provide the reader with a down-to-earth introduction and multiple access points to the circle of ideas discussed in detail in our paper \cite{EGHMW}.

Let $U$ be a connected topological space of finite homotopy type, let \begin{align*} \xi\colon\pi_1(U) \twoheadrightarrow \mathbb{Z}\end{align*} be an epimorphism, and denote by $U^\xi$ the infinite cyclic cover of $U$ corresponding to $\ker \xi$. 
Let $k$ be $\mathbb{Q}$ or $\mathbb{R}$, and denote by $R=k[t^{\pm 1}]$ the ring of Laurent polynomials in variable $t$ with $k$-coefficients.
The group $\mathbb{Z}$ of covering transformations of $U^\xi$ induces an $R$-module structure on each group $H_i(U^\xi;k)$, classically referred to as the $i$-th (homology) {\it $k$-Alexander module} of the pair $(U,\xi)$. 
Since $\xi\colon \pi_1(U) \to \mathbb{Z}$ is represented by a homotopy class of continuous maps $U \to S^1$, whenever such a representative $f\colon U \to S^1$ for $\xi$ is fixed (that is, $\xi=f_{*}$), it is also customary to use the notation $U^f$ for the corresponding infinite cyclic cover of $U$. Since $U$ is homotopy equivalent to a finite CW-complex, $H_i(U^\xi;k)$ is a finitely generated $R$-module, for each integer $i$. 

As a motivating example, let us consider the case of a fiber bundle $f\colon U \to S^1$ with connected fiber $F$ a finite CW-complex. Then $\xi=f_*\colon \pi_1(U) \to \pi_1(S^1)=\mathbb{Z}$ is surjective, and the corresponding infinite cyclic cover $U^f$ is homeomorphic to $F \times \mathbb{R}$ and hence homotopy equivalent to $F$. The deck group action on $H_i(U^f;k)$ is isomorphic 
(up to a choice of orientation on $S^1$) to the monodromy action on $H_i(F;k)$, which gives the latter vector spaces $R$-module structures. Therefore $H_i(U^f;k)\cong H_i(F;k)$ is a torsion $R$-module for all $i \geq 0$. This applies in particular to the case of the Milnor fibration $f\colon U \to S^1$ associated to a reduced complex hypersurface singularity germ, with $F$ the corresponding Milnor fiber \cite{mil}. 

However, the Alexander modules $H_i(U^\xi; k)$ are not torsion $R$-modules in general, since $U^\xi$ may not be of finite type; e.g., if $U=S^1 \vee S^2$ and $\xi=id_{\mathbb{Z}}:\pi_1(U)\cong \mathbb{Z} \to \mathbb{Z}$, then $U^\xi$ is a bouquet of $S^2$'s, one for each integer, and hence $H_2(U^\xi; k)\cong R$ . One then considers the torsion part 
\[
A_i(U^\xi;k)\coloneqq \Tors_R H_i(U^\xi;k)	
\]
of the $R$-module $H_i(U^\xi;k)$. This is a $k$-vector space of finite dimension on which a generating covering transformation (i.e., $t$-multiplication) acts as a linear automorphism.

Alexander invariants were first introduced in the classical knot theory, where it was noted that, in order to study a knot $K \subset S^3$, it is useful to consider the topology of its complement $U=S^3 \setminus K$ through the lens of the infinite cyclic cover induced by the abelianization map $\xi:\pi_1(U) \twoheadrightarrow  H_1(U; \mathbb{Z}) \cong \mathbb{Z}$. Such invariants were quickly adopted in singularity theory, for investigating the topology of Milnor fibers associated to complex hypersurface singularity germs, see, e.g., \cite{dimca1992hypersurfaces}. 

By analogy with knot theory, Libgober, Dimca, Nemethi and others considered Alexander modules for the study of the topology of complements of complex affine hypersurfaces, see, e.g., \cite{DL}, \cite{DN2},  \cite{Lib94}, \cite{Liu}, \cite{Max06}. In this context, it was shown that the Alexander modules depend on the position of singularities of a hypersurface. 

A more general setup was considered in \cite{BudurLiuWang}, where Alexander modules were associated with any complex quasi-projective manifold $U$ endowed with an epimorphism $\xi\colon \pi_1(U) \twoheadrightarrow  \mathbb{Z}$. It was shown in loc.cit. that all eigenvalues of the $t$-action on $A_i(U^\xi;k)$ are roots of unity, for any integer $i$, and upper bounds on the sizes of Jordan blocks of the $t$-action on each $A_i(U^\xi;k)$ were obtained.

In \cite{EGHMW}, we investigated a question of S. Papadima about the existence of mixed Hodge structures on the torsion parts $A_i(U^\xi;k)$ of the Alexander modules of a complex algebraic manifold $U$ endowed with an epimorphism $\xi\colon \pi_1(U) \twoheadrightarrow  \mathbb{Z}$. Note that the infinite cyclic cover $U^\xi$ is not in general a complex algebraic variety, so Deligne's classical mixed Hodge theory does not apply. 
We answered Papadima's question positively in the case when the epimorphism $\xi$ is induced by an algebraic morphism $f\colon U \to \mathbb{C}^*$. More precisely, we proved the following result (see \cite[Theorem 1.0.2]{EGHMW}):
\begin{thm}\label{mhsexistence} 
Let $U$ be a smooth connected complex algebraic variety, with an algebraic map $f\colon U \to \mathbb{C}^*$. Assume that $\xi=f_*\colon\pi_1(U) \to \mathbb{Z}$ is an epimorphism, and denote by $U^f = \{(x,z)\in U\times \C \mid f(x) = e^z \}$ the corresponding infinite cyclic cover. Then the torsion part $A_i(U^f;\Q)$ of the $\mathbb{Q}[t^{\pm 1}]$-module $H_i(U^f;\mathbb{Q})$ carries a canonical $\Q$-mixed Hodge structure for any $i\geq 0$.
\end{thm}
In this algebraic context, partial results have been previously obtained in the following special situations:
\begin{enumerate}
\item\label{hain} When $H_i(U^{\xi};\Q)$ is $\Q[t^{\pm 1}]$-torsion for all $i \geq 0$ and the $t$-action is unipotent, see \cite{hain1987rham}; 
\item \label{dlt} When $f\colon U=\mathbb{C}^n \setminus \{f=0\}\to \mathbb{C}^*$ is induced by a reduced complex polynomial $f\colon \mathbb{C}^n \to \mathbb{C}$ which is {\it transversal at infinity} (i.e., the hyperplane at infinity in $\mathbb{C}P^n$ is transversal in the stratified sense to the projectivization of $\{f=0\}$), for $\xi=f_*$ and $i<n$; see \cite{DL, Liu}. In this case it was shown in \cite{Max06} that the Alexander modules $H_i(U^\xi;\mathbb{Q})$ are torsion $\mathbb{Q}[t^{\pm 1}]$-modules for $i<n$, while $H_n(U^\xi;\mathbb{Q})$ is free and $H_i(U^\xi;\mathbb{Q})=0$ for $i>n$. Furthermore, the $t$-action on $H_i(U^\xi;\mathbb{C})$ is semisimple for $i<n$, and the corresponding eigenvalues are roots of unity of order $d=\deg(f)$.
\item\label{lkk} When $f\colon U=\mathbb{C}^n \setminus \{f=0\}\to \mathbb{C}^*$ is induced by a complex polynomial $f\colon \mathbb{C}^n \to \mathbb{C}$ which has at most isolated singularities, including at infinity, in the sense that both the projectivization of $\{f=0\}$ and its intersection with the hyperplane at infinity have at most isolated singularities. In this case, and with $\xi=f_*$, there is only one interesting Alexander module, $H_{n-1}(U^\xi;\mathbb{Q})$, which is torsion (see \cite[Theorem 4.3, Remark 4.4]{Lib94}), and a mixed Hodge structure on it was constructed in \cite{Lib96}; see also \cite{KK} for the case of plane curves under some extra conditions.
\end{enumerate}

Moreover, we showed in \cite{EGHMW} that, if $U$ and $f$ are as in case (\ref{dlt}) above, the mixed Hodge structure of Theorem \ref{mhsexistence} recovers the mixed Hodge structure obtained by different means in both \cite{DL} and  \cite{Liu}. 

In this note, we indicate the main steps in the proof of Theorem \ref{mhsexistence}, and describe several geometric applications. For complete details, the interested reader may consult our paper \cite{EGHMW}. We also rely heavily on terminology and constructions from \cite{peters2008mixed}.

\medskip

\textbf{Acknowledgements.}  E. Elduque is partially supported by an AMS-Simons Travel Grant. L. Maxim is partially supported by the Simons Foundation Collaboration Grant \#567077. B. Wang is partially supported by a Sloan Fellowship and a WARF research grant.


\section{Preliminaries}\label{s1}

\subsection{Setup. Notations. Definitions. Examples}\label{ss1}
Let $k$ be either $\mathbb{Q}$ or $\mathbb{R}$. Let $R$ denote the ring $k[t^{\pm 1}]$ of Laurent polynomials in the variable $t$ over the field $k$.

Let $U$ be a smooth connected complex algebraic variety, and let $f\colon U\rightarrow \C^*$ be an algebraic map inducing an epimorphism $f_*\colon\pi_1(U)\twoheadrightarrow \Z$ on fundamental groups. Let $\exp\colon\C\to \C^*$ be the infinite cyclic cover, and let $U^f$ be the following fiber product:
\begin{equation}\label{eq:fiberProductIntro}
\begin{tikzcd}
U^f\subset U\times \C\arrow[r,"f_\infty"] \arrow[d,"\pi"]\arrow[dr,phantom,very near start, "\lrcorner"]&
 \C \arrow[d,"\exp"] \\
 U\arrow[r,"f"] &
 \C^*.
\end{tikzcd}
\end{equation}
Hence $U^f$ is embedded in $U\times \C$ as
\[
U^f = \{ (x,z)\in U\times \C\mid f(x) = e^z\},
\]
and we let $f_\infty$ be the restriction to $U^f$ of the projection $U\times \C\to \C$. Since $\exp$ is an infinite cyclic cover, $\pi\colon U^f\to U$ is the infinite cyclic cover induced by $\ker f_*$. The group of covering transformations of $U^f$ is isomorphic to $\mathbb{Z}$, and it induces an $R$-module structure on each group $H_i(U^f;k)$, with $1\in \mathbb{Z}$ corresponding to $t \in R = k[t^{\pm 1}]$. We will also say $t$ acts on $U^f$ as the deck transformation $(x,z)\mapsto (x,z+2\pi i)$.

\begin{dfn}
The {\it $i$-th homological Alexander module} of $U$ associated to the algebraic map $f\colon U\rightarrow \C^*$ is the $R$-module
$$
H_i(U^f;k).
$$
\end{dfn}
Since $U$ has the homotopy type of a finite CW complex, the homology Alexander modules $H_*(U^f;k)$ of the pair $(U,f)$ are finitely generated $R$-modules.

\begin{ex}\label{ex1}
Let $f\colon \mathbb{C}^n \to \mathbb{C}$ be a weighted homogeneous polynomial. Then 
$f\colon U=\mathbb{C}^n \setminus \{f=0\} \to \mathbb{C}^*$ is a locally trivial topological fibration (called the global Milnor fibration of $f$), whose fiber $F$ has the homotopy type of a finite CW complex. Assume that $F$ is connected (e.g., the $\gcd$ of the exponents of the distinct irreducible factors of $f$ is $1$). As already mentioned in the introduction, in this case we have that $H_i(U^f;k)\cong H_i(F;k)$ is a torsion $R$-module for all $i \geq 0$. Moreover, the $t$-action on $H_i(U^f;k)$ is semisimple for all $i \geq 0$.
\end{ex}

\begin{ex}\label{ex2}
Let $f\colon \mathbb{C}^n \to \mathbb{C}$ be a complex polynomial with $V:=\{f=0\}$, and consider the induced map $f\colon U=\mathbb{C}^n \setminus V \to \mathbb{C}^*$. Assume $f$ has an irreducible decomposition $f=f_1^{n_1}\cdots f_r^{n_r}$ with $\gcd(n_1,\ldots, n_r)=1$. Then $H_1(U;\mathbb{Z})\cong \mathbb{Z}^r$ is generated by homology classes of positively oriented meridians $\gamma_i$ about the (regular parts of the) irreducible components of (the underlying reduced hypersurface of) $V$. Moreover, $f_*:\pi_1(U) \to \mathbb{Z}$ is the epimorphism given by assigning the integer $n_i$ to each meridian $\gamma_i$, $i=1,\ldots,r$. 

Assume moreover that $V$ is ``in general position at infinity'', that is, the hyperplane at infinity in $\mathbb{C}P^n$ is transversal (in the stratified sense) to the underlying reduced hypersurface of the projectivization of $V$. Then it was shown in \cite{MaxTommy} (see also \cite{Max06, DL, Liu} in the case when $f$ is reduced) that the Alexander modules $H_i(U^f;\mathbb{Q})$ are torsion $\mathbb{Q}[t^{\pm 1}]$-modules for $i<n$, while $H_n(U^f;\mathbb{Q})$ is free and $H_i(U^f;\mathbb{Q})=0$ for $i>n$. Furthermore, the $t$-action on $H_i(U^\xi;\mathbb{C})$ is semisimple for $i<n$, and the corresponding eigenvalues are roots of unity of order $d=\deg(f)$.
\end{ex}

\begin{ex}\label{ex1b}
Let $\cA$ be an essential hyperplane arrangement in $\C^n$, $n\geq 2$,  defined by the zero set of a reduced polynomial $f=f_1^{n_1}\cdots f_r^{n_r}:\C^n \to \C$ with $\gcd(n_1,\ldots, n_r)=1$.  Let $U=\C^n \setminus  \cA$ be the corresponding arrangement complement, and $U^f$ the infinite cyclic cover induced by $f\colon U\rightarrow \C^*$. By \cite[Theorem 4]{thesiseva}, we have that
$
H_j(U^f;\Q)$ is a torsion $R$-module for all $j<n$, a free $R$-module for $j=n$, and $0$ for $j>n$. 
\end{ex}

We next describe how the homology Alexander modules can be realized as homology groups of a certain local system on $U$. Let $ \ul k_{U^f}$ be the constant $k$-sheaf on $U^f$, and define
$$\calL = \pi_! \ul k_{U^f}.$$ The action of $t$ on $U^f$ as deck transformations induces an automorphism of $\calL$, making $\cL$ into a local system of rank $1$ free $R$-modules. It can be easily seen that $\cL$  can be described by the monodromy representation 
\[
\begin{array}{ccc}
 \pi_1(U) & \longrightarrow & \Aut_R(R)\\
 \gamma & \mapsto & \left(1\mapsto t^{f_*(\gamma)}\right)
\end{array}
\]
and, moreover, there are natural isomorphisms of $R$-modules for all $i$: 
\[
H_i(U;\calL)\cong H_i(U^f;k).
\]

For our purposes, it is more convenient to work with the cohomological version of the Alexander modules.
\begin{dfn}\label{def:2.2.6}
The {\it $i$-th cohomology Alexander module} of $U$ associated to the algebraic map $f$ is the $R$-module
$$
 H^i(U;\ov\cL),
$$
where $\ov\cL:=\cL \otimes_{t\mapsto t^{-1}} R$ is the conjugate (and dual) local system of $\cL$.
\end{dfn}

In general, the cohomological Alexander modules are not isomorphic as $R$-modules to the cohomology of the corresponding infinite cyclic cover $U^f$. Indeed, $H^i(U;\ov\cL)$ is a finitely generated $R$-module for all $i$, whereas, if $H_i(U^f;k)$ is not a finite dimensional $k$-vector space, then $H^i(U^f;k)$ is not a finitely generated $R$-module. Alternatively, if $\pi\colon U^f \to U$ is the covering map, then $H^*(U^f;k)$ and $H^*(U;\cL)$ can be computed as the cohomology of $U$ with coefficients in $R\pi_*\pi^{-1}\ul k_U=\pi_*\ul k_{U^f}$ and $\cL=\pi_!\ul k_{U^f}$ respectively, so they need not be isomorphic.


\subsection{Universal coefficient theorem and duality} \label{ss2}
The two types of Alexander modules defined in Section \ref{ss1} are related by the Universal Coefficient Theorem, in the sense that there is a natural short exact sequence of $R$-modules
\[
0\to \Ext^1_R(H_{i-1}(U;\calL),R)\to  H^i(U;\ov\cL) \to \Hom_R(H_i(U;\calL),R)\to 0.
\]
The relation between the cohomology Alexander modules and the corresponding infinite cyclic cover can be made more precise as follows.
\begin{prop}\cite[Proposition 2.4.1]{EGHMW}
\label{propcanon} \ 
There is a natural $R$-module isomorphism $$\Tors_R H^i(U;\ov\calL) \cong (\Tors_R H_{i-1}(U^f;k))^{\vee_k},$$ where $^{\vee_k}$ denotes the dual as a $k$-vector space. Moreover, this isomorphism is functorial for the pair $(U,f)$.
\end{prop}
As a consequence of Proposition \ref{propcanon}, we have the following result which applies, in particular, to the situations considered in Examples \ref{ex1},  \ref{ex2} and \ref{ex1b}.
\begin{cor} \label{isocohom}
Assume that $H_{i}(U^f;k)$ is a torsion $R$-module for some $i\geq 0$. Then, there exists a canonical isomorphism
$$
\Tors_R H^{i+1}(U;\ov\calL)\cong H^{i}(U^f;k).
$$
Moreover, if $H_{i+1}(U^f;k)$ is also a torsion $R$-module, then, so is $H^{i+1}(U;\ov\calL)$. Hence, in that case, $H^{i+1}(U;\ov\calL)$ and $H^{i}(U^f;k)$ are naturally isomorphic.
\end{cor}


\section{Mixed Hodge structures on Alexander modules}\label{s2}
In this section, we indicate the main steps in the proof of Theorem \ref{mhsexistence}. 

\subsection{Construction}\label{ss3}
A standard procedure for obtaining a mixed Hodge structure is to identify the correct mixed Hodge complex of sheaves (see, e.g., \cite[Definition 3.13]{peters2008mixed}), whose hypercohomology groups will automatically carry the desired mixed Hodge structures (see, e.g., \cite[Theorem 3.18II]{peters2008mixed}). This is roughly the approach we use in \cite{EGHMW}.

The proof of Theorem~\ref{mhsexistence} makes use of a sequence of reductions, and it relies on the construction of a suitable {\it thickening} of the Hodge-de Rham complex. We will describe these reductions and the relevant constructions in the following subsections.

We begin by noting that, since the dual of a mixed Hodge structure is again a mixed Hodge structure, the identification of Proposition \ref{propcanon} allows us to reduce the proof of Theorem \ref{mhsexistence} to the construction of mixed Hodge structures on the torsion $R$-modules $\Tors_R H^i(U;\ov\calL)$.

\subsubsection{Reduction to the case of unipotent $t$-action}\label{sss4}
Let $k=\C$, let $U$ be a smooth connected complex algebraic variety, and let $f\colon U\rightarrow \C^*$ be an algebraic map with associated local system $\calL$ of $R$-modules as in Section~\ref{ss1}. Since $R$ is a PID, we have the primary decomposition
$$
A_i(U^f;\C):=\Tors_R H_i(U;\calL)\cong \bigoplus_{j=1}^r R/\big((t-\lambda_j)^{p_j}\big)
$$
with $p_j\geq 1$ for all $j=1,\ldots, r$. The set $\{\lambda_j \in \C \mid j=1, \ldots, r\}$ is uniquely determined by $A_i(U^f;\C)$.
The following result was essentially proved in \cite[Proposition 1.4]{BudurLiuWang}, but see also \cite[Proposition 2.6.1]{EGHMW} for a slight generalization to the current setup.
\begin{prop}\label{eig1}
Every $\lambda_j$ defined above is a root of unity.
\label{monodromy}
\end{prop}
In particular, using Proposition \ref{propcanon}, one has the following.
\begin{cor}\label{abc}
Let $k=\C$. The eigenvalues of the action of $t$ on $\Tors_R H^*(U;\ov\calL)$ are all roots of unity.
\label{roots}
\end{cor}

Corollary \ref{abc} implies that we can choose $N\in\N$ such that $t^N-1$ acts nilpotently on $\Tors_R H^i(U;\ov\calL)$ for all $i$. Consider the following pull-back diagram:
$$
\begin{tikzcd}
 U_N = \{(x,z)\in U\times \C^*\mid f(x) = z^N \}\arrow[d,"p"] \arrow[dr,phantom,very near start,"\lrcorner"] \arrow[r,"f_N"] & \C^*\arrow[d,"z\mapsto z^N"] \\
 U \arrow[r,"f"] & \C^*
\end{tikzcd}
$$
where $p$ is an $N$-sheeted cyclic cover, and note that all maps involved in this diagram are algebraic and $U_N$ is a smooth algebraic variety. We can then define, as in Section~\ref{ss1},  $(U_N)^{f_N}$, $(f_N)_\infty$, $\pi_N$ and $\cL_N$  for the map $f_N\colon U_N\to \C^*$. We also define:
\[
\f{\theta_N}{U^f}{U_N^{f_N}}{U\times \C\ni (x,z)}{(x,e^{z/N}, z/N)\in U_N\times \C\subset U\times \C^* \times \C,}
\]
which fits into the following commutative diagram:
\begin{equation}\label{eq:UN}
\begin{tikzcd}[column sep = 7em]
U^f \arrow[d,"f_\infty"]\arrow[r,"\sim"',"\theta_N", dashed] 
\arrow[rrr,"\pi",rounded corners,to path=
{ --([yshift = 0.7em]\tikztostart.north)
-- ([yshift = 0.7em]\tikztotarget.north)\tikztonodes
-- (\tikztotarget)}] 
&
U^{f_N}_N \arrow[d,"(f_N)_\infty"]\arrow[r,"\pi_N"]\arrow[dr,phantom,very near start,"\lrcorner"] &
U_N \arrow[d,"f_N"] \arrow[r,"p"]\arrow[dr,phantom,very near start,"\lrcorner"] &
U \arrow[d,"f"] \\
\C \arrow[r,"z\mapsto\frac{z}{N}"]\arrow[rrr,"\exp", rounded corners,to path=
{ --([yshift = -0.7em]\tikztostart.south)
-- ([yshift = -0.7em]\tikztotarget.south)\tikztonodes
-- (\tikztotarget)}] &
\C \arrow[r,"\exp"]&
\C^* \arrow[r,"z\mapsto z^N"]&
\C^*.
\end{tikzcd}
\end{equation}
The map $\theta_N$ allows us to identify $U^f$ with $U_N^{f_N}$ in a canonical way, which we will do from now on. In particular, we can also identify the constant sheaves $\ul k_{U_N^{f_N}}$ and $\ul k_{U^f}$ canonically.

Let $R(N)\coloneqq k[t^{N},t^{-N}]$. Since the deck group of the infinite cyclic cover $\pi_N$ is generated by $t^N$, the corresponding Alexander modules $H_i(U_N^{f_N};k)$ are finitely generated $R(N)$-modules. Since $R$ is a rank $N$ free $R(N)$-module, we can also consider $\calL$ as a local system of rank $N$ free $R(N)$-modules on $U$. Moreover, $\theta_N$ induces the canonical isomorphism $p_*\cL_N = \cL$ of local systems of $R(N)$-modules, which can be further used to prove the following (cf. \cite[Proposition 2.6.3]{EGHMW}).
\begin{lem}\label{233}
\label{lemLocal}
In the above notations, $\theta_N$ induces the following canonical isomorphisms of $R(N)$-modules:
\[\Tors_R H^{i}(U;\ov\cL) \cong \Tors_{R(N)}H^{i}(U_N;\ov\cL_N) \]
for any integer $i \geq 0$.
\end{lem}

\begin{remk}
Notice that the only eigenvalue of the action of $t^N$ in $\Tors_{R(N)} H^*(U_N;\ov\calL_N)$ is $1$. So Lemma \ref{233} allows us to reduce the problem of constructing a mixed Hodge structure on $\Tors_{R} H^*(U;\ov\calL)$ to the case when the only eigenvalue is $1$.
\label{remEigenvalue}
\end{remk}


\subsubsection{Isolating the torsion}\label{sss5}
As seen in Section \ref{sss4}, upon passing to a finite cover $U_N$ of $U$ (and replacing $f$ by $f_N$, $\cL$ by $\cL_N$ and $t$ by $t^N$), we may assume that $t-1$ acts nilpotently on  $\Tors_{R} H^i(U;\ov\calL)$, for any integer $i \geq 0$. We will assume this is the case from now on.

In this section, we explain how  $\Tors_{R} H^i(U;\ov\calL)$ can be understood in terms of the cohomology of certain $k$-local systems of finite rank on $U$. 

Let $s=t-1$. For $m \in \mathbb{N}$, we set \[ R_m:=R/s^mR \] and let \[ \ov\calL_m:=\ov\calL \otimes_R R_m \] be the corresponding rank one local system of $R_m$-modules.
Since $s$ acts nilpotenly on $\Tors_{R} H^*(U;\ov\calL)$, there exists an integer $m\geq 0$ such that $s^m$ annihilates $\Tors_{R} H^i(U;\ov\calL)$, for all $i \geq 0$. With the above notations, it can be seen (cf. \cite[Lemma 3.1.8,  Corollary 3.1.9]{EGHMW}) that the maps of sheaves
\[
\ov\calL \twoheadrightarrow \ov\calL_m \overset{\cdot s^m}{\hookrightarrow} \ov\calL_{2m}
\]
induce an exact sequence
\[
0 \to \Tors_{R} H^*(U;\ov\calL)  \to H^*(U;\ov\calL_m) \overset{\cdot s^m}{\to}  H^*(U;\ov\calL_{2m}).
\]
Hence, 
\begin{equation}\label{ker} 
\Tors_{R} H^*(U;\ov\calL)  \cong \ker \left( H^*(U;\ov\calL_m) \overset{\cdot s^m}{\to}  H^*(U;\ov\calL_{2m}) \right)  
\end{equation}

Our next goal is to endow each $H^i(U;\ov\calL_m)$ with a canonical mixed Hodge structure for all $m \geq 1$, such that the map 
\[ H^*(U;\ov\calL_m) \overset{\cdot s^m}{\to}  H^*(U;\ov\calL_{2m})(-m) \]
 is a morphism of mixed Hodge structures (where $(-m)$ denotes the $-m$th Tate twist). This will be achieved by resolving $\ov\calL_m$ by a certain mixed Hodge complex. Note that a mixed Hodge structure on $H^*(U;\ov\calL_m)$  induces by \eqref{ker} a canonical mixed Hodge structure on $\Tors_{R} H^*(U;\ov\calL)$.


\subsubsection{Thickening Deligne's Hodge-de Rham mixed Hodge complex}
Assume $k = \mathbb{R}$, and let as before $s=t-1$. 
Let $(X,D)$ be a good compactification of the smooth variety $U$, and let $j:U \hookrightarrow X$ be the inclusion map. Let $\calE^\bullet_U$ be the real de Rham complex, and  $\logdr{X}{D}$ the log de Rham complex (see, e.g., \cite[Section 4.1]{peters2008mixed}). It is known by work of Deligne that $(j_*\calE^\bullet_U, \logdr{X}{D})$ form the real and, resp., complex part of the {\it Hodge-de Rham mixed Hodge complex}
\[
\calH dg^\bullet(X\,\log D),
\]
whose hypercohomology endows $H^*(U;\mathbb{R})$ with a canonical mixed Hodge structure.

To resolve $\ov\calL_m:=\ov\calL \otimes_R R_m$ by a mixed Hodge complex, we will perform a thickening of  the Hodge-de Rham mixed Hodge complex, a procedure already used in \cite{BudurLiuWang}. 
Recall that for any $m \geq 1$, an {\it $m$-thickening} of a $k$-cdga $(A,d,\wedge)$ in the direction $\eta \in A^1 \cap \ker d$ is the cochain complex of $R_m$-modules denoted by 
\begin{align*}
	A(\eta,m) = (A {\otimes_k} R_m, d_\eta)
\end{align*}
 and described by:
\begin{itemize}
\item[(i)] for $p \in \Z$, the $p$-th graded component of $A(\eta,m)$ is {$A^p {\otimes_k} R_m$.}
\item[(ii)] for $\omega \in A$ and $\phi \in R_m$, we set $d_\eta(\omega \otimes \phi) = d\omega \otimes \phi + (\eta \wedge \omega) \otimes {s\phi} $.
\end{itemize}
If $(\calA, \wedge, d)$ is a sheaf of cdgas on $X$ and $\eta \in \Gamma(X, \calA^1) \cap \ker d$ is a closed global section, then we define similarly 
 the {$m$-thickening of $\calA$ in direction $\eta$}, denoted by $\calA(\eta, m)$.

In the above notations, we have the following result (see \cite[Section 5.2]{EGHMW}).
\begin{lem}\label{ltic} Let $k = \mathbb{R}$ and let $\calE^\bullet_U$ be the real de Rham complex on $U$. Then a canonical resolution of $\ov\calL_m$ as a sheaf of $R_m$-modules is given by the thickened complex $\calE^\bullet_U(\Im \frac{df}{f},m)$ with a modified $R_m$-module structure, the action of $s$ becoming multiplication by $-\frac{\log t}{2\pi}$ (expressed as a power series in $s=t-1$). Here, $\Im$ denotes the imaginary part.
\end{lem}

We can similarly thicken the log de Rham complex $\logdr{X}{D}$ by the logarithmic form $\frac{1}{i}\frac{df}{f}$ (which is cohomologous to $\Im \frac{df}{f}$). Puting everything together, this leads to a thickening \[ \calH dg^\bullet(X\,\log D)\left(\frac{1}{i}\frac{df}{f},m\right) \] of the Hodge-de Rham complex, whose hypercohomology computes $H^*(U;\ov\calL_m)$. Moreover, since by \cite[Theorem 4.2.1]{EGHMW} the thickened complex of a mixed Hodge complex of sheaves is again a mixed Hodge complex of sheaves (assuming $\eta \in W_1 \cap F^1$), the thickened  Hodge-de Rham complex $ \calH dg^\bullet(X\,\log D)\left(\frac{1}{i}\frac{df}{f},m\right)$ is an $\R$-mixed Hodge complex of sheaves on $X$ (cf. \cite[Theorem 5.4.3]{EGHMW}).
This then yields:
\begin{cor} For all $i$ and $m$, 
\[ H^i(U;\ov\calL_m) \cong \mathbb{H}^i\left(X, j_*\calE^\bullet_U\left(\Im \frac{df}{f},m\right)\right) \]
has a canonical $\mathbb{R}$-mixed Hodge structure. 
\end{cor}
While we omit here the discussion on filtrations (roughly speaking, these are made of the filtrations on the Hodge-de Rham complex and powers of $s$), let us just say that these are defined so that multiplication by $s$ induces a mixed Hodge complex morphism into the $-1$st Tate twist. In view of \eqref{ker} and Lemma \ref{ltic} this then yields:
\begin{cor}\label{torsionmhs}
Suppose that the action of $t$ on $\Tors_{R}H^*(U; \ov\calL)$ is unipotent. The $\R$-vector spaces $\Tors_{R}H^*(U; \ov\calL)$ admit natural $\R$-mixed Hodge structures for which multiplication by $\log t$ (seen as a power series in $s=t-1$) determines a morphism of mixed Hodge structures into the $-1$st Tate twist.
\end{cor}
\begin{remk}\label{logtN}
If the $t$-action on  on $\mathrm{Tors}_R\, H^*(U;\ov\calL)$ is not unipotent, then $t$ gets replaced by $t^N$ which acts unipotently, and then $\log t^N$  is a morphism of mixed Hodge structures into the $-1$st Tate twist. 
\end{remk}

Finally, we have the following result (cf. \cite[Theorem 5.4.10]{EGHMW}) which, in view of Proposition \ref{propcanon}, completes the justification of Theorem 
\ref{mhsexistence}:
\begin{thm}[$\Q$-MHS]\label{Qalexandermhs}
The mixed Hodge structure on $\Tors_R H^*(U;\ov\cL)$ defined above for $k=\R$ comes from a (necessarily unique) mixed Hodge structure defined for $k=\Q$.
\end{thm}


\subsection{Properties}

Let $U$ be a smooth connected complex algebraic variety, and let $f\colon U\rightarrow \C^*$ be an algebraic map inducing an epimorphism $\pi_1(U) \twoheadrightarrow \Z$, with corresponding infinite cyclic cover $U^f$. There are several choices made in the construction of our mixed Hodge structure on the torsion part $A_*(U^f;\Q)$ of the Alexander modules, e.g., a good compactification $(X,D)$ of $U$, a finite cover $U_N$ of $U$ on which the monodromy is unipotent, a positive integer $m$ such that $(t-1)^m$ annihilates $\Tors_R H^*(U;\ov\cL)$. In \cite{EGHMW}, we showed that the mixed Hodge structure we constructed on  $A_*(U^f;\Q)$ is independent of all these choices, see  \cite[Theorem 5.4.7]{EGHMW} (independence of compactification),  \cite[Theorem 5.4.8]{EGHMW} (independence of $U_N$), \cite[Corollary 5.4.4]{EGHMW} (independence of nilpotence index $m$). Furthermore, the mixed Hodge structure on $A_*(U^f;\Q)$
behaves functorially with respect to algebraic maps over $\C^*$, see \cite[Theorem 5.4.9]{EGHMW}. 

In view of Proposition \ref{propcanon}, it suffices to prove the above statements for the mixed Hodge structure on $\Tors_R H^*(U;\ov\cL)$, in which case the assertions follow, roughly speaking, by constructing appropriate morphisms of mixed Hodge complexes of sheaves and taking hypercohomology. 


\section{Relation with other mixed Hodge structures and applications}
In this section, we relate our mixed Hodge structures on the torsion parts of the Alexander modules with other known mixed Hodge structures, and indicate several geometric consequences. 

\subsection{Hodge theory of the infinite cyclic covering map. Applications}
The infinite cyclic covering map $U^f \to U$ induces a natural map of vector spaces $A_i(U^f; \Q) \to H_i(U;\Q)$. In \cite[Theorem 6.0.1]{EGHMW}, we prove the following result.

\begin{thm}\label{geoIntro}
In the setting of Theorem~\ref{mhsexistence}, the vector space homomorphism \[ H_i(\pi)\colon A_i(U^f; \Q) \to H_i(U;\Q) \] induced by the covering $\pi\colon U^f \rightarrow U$ is a morphism of mixed Hodge structures for all $i \geq 0$, where $H_i(U;\Q)$ is equipped with (the dual of) Deligne's mixed Hodge structure. 
\end{thm}

Using the faithful flatness of $\R$ over $\Q$ it suffices to work with $k=\R$ and prove the result in cohomology. Once again, this amounts to constructing an appropriate morphism of mixed Hodge complexes of sheaves which, upon taking hypercohomology, gives a mixed Hodge structure morphism $H^i(\pi)\colon H^i(U;\R) \to \Tors_R H^{i+1}(U;\ov\cL)$. The geometric meaning of the morphism $H^i(\pi)$ is achieved by carefully tracing arrows in the derived category. 

\begin{remk}
The notation $H^i(\pi)$ of the previous paragraph is justified by the fact that, under the assumption that $H_i(U^f;\R)$ is a torsion $R$-module, the above $H^i(\pi)$ coincides with the map $H^i(U;\R) \to H^i(U^f;\R)$ induced in cohomology by the covering map $\pi\colon U^f \rightarrow U$ (the dual of the map $H_i(\pi)$ induced in homology); compare with Corollary \ref{isocohom}.
\end{remk}

In the following subsection we indicate several applications of Theorem \ref{geoIntro}.

\subsubsection{Bounding the weights. Size of Jordan blocks}
 Theorem~\ref{geoIntro} and our construction can be used to obtain a bound on the weight filtrations of the mixed Hodge structures on the torsion parts of the Alexander modules (see \cite[Theorem 7.4.1]{EGHMW}). This bound coincides with the known bound for the homology of smooth algebraic varieties of the same dimension as the generic fiber of $f$ (cf. \cite[Corollaire 3.2.15]{De2}). Specifically, we get the following.

\begin{thm}\label{boundsIntro}
Assume the setting of Theorem~\ref{mhsexistence}. Let $i \geq 0$.
If $\ell \notin [i,2i]\cap [i,2\dim_\C(U)-2]$, then
\[
\Gr^W_{-\ell} A_i(U^f;\Q) = 0
\]
where $\Gr^W_{-\ell}$ denotes the $-\ell$th graded piece of the weight filtration.
\end{thm}

This theorem is proved by first using Lemma~\ref{233} to reduce to the case where $t$ acts unipotently on $A_i(U^f;\Q)$ for all $i$. The main step in the proof amounts to constructing the following exact sequence of mixed Hodge structures:
\[
H^i(U;\Q)\xrightarrow{H^i(\pi)} \Tors_R H^{i+1}(U;\ov\cL)\xrightarrow{\cdot \log (t)}  \Tors_R H^{i+1}(U;\ov\cL)(-1)\to H^{i+1}(U;\Q).
\]
Given this exact sequence, the assertion of Theorem \ref{boundsIntro} follows by using the bounds for the nonzero weights on the cohomology of the smooth variety $U$ together with the finite dimensionality of the mixed Hodge structures.

\medskip

Other consequences of our construction and Theorem~\ref{geoIntro} are related to the $t$-action on the torsion parts of the Alexander modules. 
For example, we apply it to determine bounds on the size of the Jordan blocks of this $t$-action (see \cite[Corollary 7.4.2]{EGHMW}), which nearly cut in half existing bounds, as in \cite[Proposition 1.10]{BudurLiuWang}. Specifically, we prove:
\begin{cor}\label{cor:jordan}
Every Jordan block of the action of $t$ on $A_i(U^f;\C)$ 
has size at most \[ \min\{ \lceil (i+1)/2\rceil, n- \lfloor (i+1)/2\rfloor\}.\] In particular, $A_1(U^f;\C)$ is a semisimple $R$-module.
\end{cor}

To prove the above corollary, one simply needs to observe that, by Remark~\ref{logtN}, applying $\log(t^N)$ decreases the weight by 2, and use Theorem~\ref{boundsIntro}. This implies that $(\log(t^N))^m=0$, where $m$ is the bound in the corollary. Finally one observes that this implies that $(t^N-1)^m=0$.

\begin{remk}\label{remk:A1Semisimple}
If $f\colon \C^2\rightarrow \C$ is a polynomial function such that $f^{-1}(0)$ is reduced and connected, let $U=\C^2\setminus f^{-1}(0)$, and consider the induced map $f\colon U\rightarrow \C^*$. By \cite[Corollary 1.7]{DL}, the action of $t$ on $H_1(U^f;\C)$ is semisimple. The last part of the above corollary generalizes this result, not only to algebraic varieties $U$ that are not affine connected curve complements, but also to connected curve complements in which the corresponding $f$ is given by a non-reduced polynomial, or even a rational function.
\end{remk}


\subsubsection{Semisimplicity and its consequences} A very natural question is to understand under which conditions the $t$-action on $A_i(U^f; \Q)$  is a mixed Hodge structure morphism. We prove the following result (see \cite[Corollary 7.0.4 and Proposition 7.0.5]{EGHMW}):

\begin{thm}\label{tsemisimpleIntro}
	Assume the setting of Theorem \ref{mhsexistence}. Let $i \geq 0$. The $t$-action on $A_i(U^f; \Q)$ is a mixed Hodge structure morphism if and only if it is semisimple.
\end{thm}

The `only if' portion of the above theorem is a consequence of a simple fact about mixed Hodge structures: if a mixed Hodge structure $V$ is endowed with an endormophism $t$ such that $\log(t^N)$ is a morphism $V\to V(-1)$ (Remark~\ref{logtN}), then necessarily $t^N=\id_V$. For the forward direction, it suffices to use note that $t^N$ is unipotent for some $N$ (Proposition~\ref{eig1}), so the hypothesis implies that $t^N=\id$. The Milnor long exact sequence then shows that $A_i(U^f; \Q)$ is a quotient of $H_i(U_N;\Q)$. It only remains to note that $t$ acts on the latter group by deck transformations on $U_N$, an algebraic map, which makes the induced map in homology a mixed Hodge structure morphism.

\medskip

When the $t$-action is semisimple, we show that the mixed Hodge structure on the torsion parts $A_i(U^f;\Q)$ of the Alexander modules can be constructed directly using a finite cyclic cover, which, unlike an infinite cyclic cover, is always a complex algebraic variety.
This bypasses our rather abstract general construction of the mixed Hodge structure. 
In \cite{EGHMW}, we present two different viewpoints.
In the first, we utilize cap product with the pullback of a generator of $H^1(\C^*;\Q)$. 
In the second, we utilize a generic fiber of the algebraic map, which is always a complex algebraic variety. For the following result see \cite[Corollary 7.1.3 and Corollary 7.2.1]{EGHMW}.

\begin{thm}\label{finiteIntro}
	Assume the setting of Theorem \ref{mhsexistence}. 
	Let $i \geq 0$ and assume that the $t$-action on $A_i(U; \Q)$ is semisimple. 
	Let $N$ be such that the action of $t^N$ on $A_i(U^f;\Q)$ is the identity, and let $U_N = \{(x,z) \in U \times \C^*\mid f(x) = z^N\}$ denote the corresponding $N$-fold cyclic cover.
	Equip the rational homology of $U_N$ with the (dual of) Deligne's mixed Hodge structure.
	
	\begin{enumerate}[(A)]
	\item Let $f_N\colon U_N \rightarrow \C^*$ denote the algebraic map induced by projection onto the second component, and let $\gen \in H^1(\C^*;\Q)$ be a generator. Then $A_i(U^f;\Q)$ is isomorphic as a mixed Hodge structure to the image of the mixed Hodge structure morphism induced by cap product with $f^*_N(\gen)$
	\[(-) \frown f^*_N(\gen) \colon H_{i+1}(U_N;\Q)(-1) \rightarrow H_i(U_N;\Q) ,\]
	where $(-1)$ denotes the $-1$th Tate twist of a mixed Hodge structure. 
	
	\item Let $F \hookrightarrow U$ be the inclusion of any generic fiber of $f$ and let $F \hookrightarrow U_N$ be any lift of this inclusion. 
	Then $A_i(U^f;\Q)$ is isomorphic as a mixed Hodge structure to the image of the mixed Hodge structure morphism
	\[H_i(F; \Q) \rightarrow H_i(U_N; \Q)\] 
	induced by the inclusion, where $H_i(F;\Q)$ is equipped with (the dual of) Deligne's mixed Hodge structure.
	\end{enumerate}
\end{thm}

The first viewpoint granted by semisimplicity, in terms of cap products (Theorem \ref{finiteIntro}A), is suggested by the thickened complexes that play the central role in our construction. 

As for the second viewpoint, note that the homologies of different choices of generic fibers in the same degree may have different mixed Hodge structures, but \emph{any} choice is allowed in Theorem \ref{finiteIntro}B.
This shows that the mixed Hodge structures on the torsion parts of the Alexander modules are common quotients of the homologies of all generic fibers, when semisimplicity holds. In fact, under the assumptions of Theorem \ref{finiteIntro}, the inclusion $F \hookrightarrow U $ of a generic fiber of $f\colon U \to \C^*$ lifts to $U_N$ and $U^f$ via maps $i_N$ and $i_\infty$, making the following diagram commutative, where the vertical arrows are covering space maps.
$$
\begin{tikzcd}
\ & U^f\arrow[d, "\pi_N"]\arrow[dd, "\pi", bend left=90]\\
\ & U_N\arrow[d,"p"]\\
F\arrow[r,"i",hook]\arrow[ru, "i_N" pos=1, hook]\arrow[ruu, "i_{\infty}",hook] & U
\end{tikzcd}
$$
Note that, in homology, the composition $i_N=\pi_N\circ i_{\infty}$ factors through $A_*(U^f;\Q)$, hence we get a diagram:
\begin{equation}\label{eq:fiberhom}
\begin{tikzcd}
H_j(F;\Q)\arrow[r,"H_j(i_\infty)"]\arrow[rr, "H_j(i_N)", bend right = 20] & A_j(U^f;\Q)\arrow[r, "H_j(\pi_N)"] & H_j(U_N;\Q).
\end{tikzcd}
\end{equation}
Moreover, it can be shown by topological arguments that  $H_j(i_\infty)$ is surjective (see \cite[Proposition 2.5.3]{EGHMW}), and the semisimplicity assumption on the $t$-action yields that $H_j(\pi_N)$ is injective (see \cite[Corollary 7.0.2]{EGHMW}). Since $i_N$ is an algebraic map, it follows that $H_j(i_N)$ is a morphism of mixed Hodge structures, and the same is true for $H_j(\pi_N)$ by Theorem \ref{geoIntro}. Putting everything together, we get the following more general consequence of Theorem \ref{geoIntro} (see \cite[Corollary 7.2.1]{EGHMW}), from which Theorem \ref{finiteIntro}B follows readily.
\begin{cor}\label{cor:fiber}
Let $N$ be such that the action of $t^N$ on $A_j(U^f;\Q)$ is unipotent. Suppose that the $t$-action on $A_j(U^f;\Q)$ is semisimple. Then, we have the following commutative diagram, where all the arrows are morphisms of mixed Hodge structures:
$$
\begin{tikzcd}
H_j(F;\Q)\arrow[r,"H_j(i_\infty)",two heads]\arrow[rr, "H_j(i_N)", bend right = 20] & A_j(U^f;\Q)\arrow[r, "H_j(\pi_N)", hook] & H_j(U_N;\Q).
\end{tikzcd}
$$
\end{cor}

Making use of Theorem~\ref{tsemisimpleIntro}, we also have the following (see \cite[Corollary 7.2.3]{EGHMW}).
\begin{cor}\label{cor:fiberConverse}
The $t$-action on $A_j(U^f;\Q)$ is semisimple if and only if for any generic fiber $F\subset U^f$, the induced map in homology $H_j(i_\infty)\colon H_j(F;\Q)\to A_j(U^f;\Q)$ is a mixed Hodge structure morphism.
\end{cor}

\begin{ex}[Global Milnor fiber]\label{ex4}
Let $f\in\C[x_1,\ldots,x_n]$ be a weighted homogeneous polynomial, and let $U=\C^n\setminus\{f=0\}$. As already mentioned in Example \ref{ex1}, we have a global Milnor fibration $f\colon U\to\C^*$. Assume that $f\colon U\to\C^*$ induces an epimorphism on fundamental groups, i.e., the greatest common divisor of the exponents of the distinct irreducible factors of $f$ is $1$. Let $F$ be a fiber of $f\colon U\to\C^*$, and let as before $i_\infty\colon F\hookrightarrow U^f$ be a lift of the inclusion $i\colon F\hookrightarrow U$. Since $f\colon U\rightarrow \C^*$ is a fibration, we have that $i_\infty$ is a homotopy equivalence, so it induces isomorphisms $H_j(F;\Q)\rightarrow H_j(U^f;\Q)$ for all $j$, which are compatible with the $t$-action (see \cite[Lemma 2.5.2]{EGHMW}). Since the $t$-action on $F$ comes from an algebraic map $F\rightarrow F$ of finite order, it follows that the $t$-action on $H_j(F)$ is semisimple. Applying Corollary \ref{cor:fiber}, we see that the Alexander modules recover in this case the Deligne mixed Hodge structure on the global Milnor fiber. Specifically, the map
\[ H_j(F;\Q)\rightarrow H_j(U^f;\Q) \] induced by $i_\infty$ is a mixed Hodge structure morphism, where $H_j(U^f;\Q)$ is endowed with the (dual) Deligne mixed Hodge structure.
\end{ex}

\begin{remk}
Recall from Example \ref{ex2} that if $f:\C^n \to \C$ is a reduced complex polynomial so that $V=\{f=0\}$ is in general position at infinity, then 
the Alexander modules $H_i(U^f;\mathbb{Q})$ of $U=\C^n \setminus V$ and $f\colon U \to \C^*$ are torsion $\mathbb{Q}[t^{\pm 1}]$-modules for $i<n$, while $H_n(U^f;\mathbb{Q})$ is free and $H_i(U^f;\mathbb{Q})=0$ for $i>n$. Moreover,  the $t$-action on $H_i(U^\xi;\mathbb{C})$ is semisimple for $i<n$, and the corresponding eigenvalues are roots of unity of order $d=\deg(f)$. A mixed Hodge structure on $H_i(U^\xi;\mathbb{Q})$,  for $i<n$, has been constructed by Dimca-Libgober \cite{DL} and Liu \cite{Liu}. 
Corollary \ref{cor:fiber} can be used directly to show that our mixed Hodge structure on $A_*(U^f;\Q)$ coincides in this case with those constructed by Dimca-Libgober and Liu, see \cite[Corollary 7.3.6]{EGHMW}.
\end{remk}

Theorem \ref{finiteIntro}, when it applies, reinforces the significance of semisimplicity. Our results from \cite{EGHMW} show that, in fact, semisimplicity is not a rare occurrence (we have already encountered such instances in Examples \ref{ex1} and \ref{ex2}). For instance, when $f$ is proper, we have the following (see \cite[Corollary 8.0.2]{EGHMW}):
\begin{thm}\label{thmsimple}
Let $U$ be a smooth complex algebraic variety, and let $f\colon U\to \C^*$ be a proper algebraic map. Then the torsion part $A_i(U^f; \Q)$ of the homology Alexander module $H_i(U^f; \Q)$ is a semisimple $R$-module, for all $i \geq 0$. 
\end{thm}

If the map $f\colon U \to \C^*$ of Theorem \ref{thmsimple} is a projective submersion, then $f$ is a fibration, and let $F$ be its fiber. The semisimplicity of $A_i(U^f;\Q) \cong H_i(U^f;\Q) \cong H_i(F;\Q)$ is in this case a direct consequence of Deligne's decomposition theorem \cite{De, De2}. In the general case, Theorem \ref{thmsimple} is proved by using the decomposition theorem of Be\u{\i}linson--Bernstein--Deligne \cite{BBD}.

In view of Corollary \ref{cor:fiber}, a nice application of the semisimplicity statement of Theorem \ref{thmsimple} is the following purity result (see \cite[Corollary 8.0.6]{EGHMW}).
\begin{thm}\label{thm:pur} If $f\colon U \to \C^*$ is a proper algebraic map, then $A_i(U^f;\Q)$ carries a pure Hodge structure of weight $-i$.
\end{thm}

\begin{remk}
In fact, we do not know of any example where semisimplicity does not hold. This lack of examples is mainly due to the fact that higher Alexander modules are harder to compute than the first (which, as seen in Corollary \ref{cor:jordan}, is always semisismple, and can be computed from a presentation of the fundamental group). 
\end{remk}


\subsection{Relation with the limit mixed Hodge structure}
The mixed Hodge structure on the torsion part $A_*(U^f;\Q)$ of the Alexander modules of the pair $(U, f)$ can be regarded as a global version of the {\it limit mixed Hodge structure} on the generic fiber of $f$, in the following sense. 
Let $f\colon U \to \C^*$ be an algebraic map inducing an epimorphism on fundamental groups, and let $U^f$ denote as before the corresponding infinite cyclic cover of $U$. Let $D^*$ be a sufficiently small punctured disk centered at $0$ in $\C$, such that $f\colon f^{-1}(D^*)\rightarrow D^*$ is a fibration, and let $T^*=f^{-1}(D^*)$. The infinite cyclic cover $(T^*)^f$ is homotopy equivalent to $F$, where $F$ denotes any fiber of the form $f^{-1}(c)$, for $c\in D^*$. In fact, $(T^*)^f$ can be regarded as the canonical fiber of $f\colon f^{-1}(D^*)\rightarrow D^*$. If $f$ is proper, $H_i((T^*)^f;\Q)$ is also endowed with a limit mixed Hodge structure, which can be compared with the one we constructed on $A_i(U^f;\Q)$ via the following result:

\begin{thm}\label{comp} In the setup of Theorem \ref{mhsexistence}, assume moreover that 
$f\colon U \to \C^*$ is proper. Then, in the above notations,  the inclusion $(T^*)^f\subset U^f$ induces for all $i\geq 0$ an epimorphism of $\Q$-mixed Hodge structures
\begin{align}\label{mhscomp}
H_i((T^*)^f;\Q) \twoheadrightarrow A_i(U^f;\Q),	
\end{align}
where $H_i((T^*)^f;\Q)$ is endowed with its limit mixed Hodge structure. If, moreover, $f$ is a fibration, then the two mixed Hodge structures are isomorphic.
\end{thm}
The mixed Hodge structure morphism of Theorem \ref{comp} is realized, upon taking hypercohomology and $\Q$-duals, by a suitable morphism of mixed Hodge complexes of sheaves. In more detail, let $(X, D)$ be a good compactification  of $U$ by a simple normal crossing divisor $D = X\setminus U$ such that $f\colon U \rightarrow \C^*$ extends to an algebraic map $\bar{f}\colon X \rightarrow \C P^1$.
By replacing $f\colon U \rightarrow \C^*$ with a finite cyclic cover $f_N\colon U_N \rightarrow \C^*$ if necessary, we may assume that $E:=\bar{f}^{-1}(0)$ is reduced and $X \setminus U = \bar{f}^{-1}(\{0, \infty\})$. Let $i\colon E \hookrightarrow X$ be the inclusion. By restricting $\bar{f}$ above a sufficiently small punctured disk $D^*$ centered at $0 \in \C$, one can define the nearby cycle functor $\psi_{\bar{f}}$ of Deligne, and there is a vector space isomorphism  
\[ \mathbb{H}^*(E; \psi_{\bar{f}}\underline{\Q}) \cong H^*(F; \Q),\] where $F$ is any fiber of $f$ over $D^*$. A clockwise loop in $D^*$ determines a monodromy homeomorphism from $F$ to itself and so equips $\mathbb{H}^*(E; \psi_{\bar{f}}\underline{\Q})$ with the structure of a torsion module over $\Q[t^{\pm 1}]$. The limit mixed Hodge  structure on $H^*(F;\Q)$ is realized by a mixed Hodge complex $\psi_{\bar{f}}^{\textnormal{Hdg}} $ of sheaves on $E$, assigned to $\psi_{\bar{f}}\underline{\Q}$; see \cite[Theorem 11.22]{peters2008mixed}, and also \cite[Theorem 2.11.1]{EGHMW}. The mixed Hodge structure morphism \eqref{mhscomp} is then induced by a morphism of mixed Hodge complexes from the thickened Hodge-de Rham complex (shifted by $[1]$ and with an appropriate twisting of the $R$-module structure) and $i_*\psi_{\bar{f}}^{\textnormal{Hdg}}$. For complete details and a geometric interpretation of this morphism of mixed Hodge complexes, see \cite[Section 9]{EGHMW}.

\section{Examples. Hyperplane arrangements}
Let $n\geq 2$. Let $f_1,\ldots,f_d$ be degree $1$ polynomials in $\C[x_1,\ldots,x_n]$ defining $d$ distinct hyperplanes and let $f=f_1\cdot\ldots\cdot f_d$. The zeros of $f$ define a hyperplane arrangement $\cA$ of $d$ hyperplanes in $\C^n$. Let $U\subset \C^n$ be the corresponding arrangement complement, with induced map $f\colon U \to \C^*$. For the purpose of studying Alexander invariants of the pair $(U,f)$, it suffices to assume that $\cA$ is essential, that is, the intersection of some subset of hyperplanes of $\cA$ is a point (see, e.g., \cite[Remark 10.1.2]{EGHMW}).  
As already indicated in Example \ref{ex1b}, if $\cA$ is essential then $H_j(U^f;\Q)$ is a torsion $R$-module for all $j<n$, a free $R$-module for $j=n$, and $0$ for $j>n$. In particular, by Theorem  \ref{mhsexistence}, we can endow $H_j(U^f;\Q)$  and $H^j(U^f;\Q)$ with canonical mixed Hodge structures, for $0\leq j \leq n-1$.

By Corollary \ref{cor:jordan}, the $t$-action on $H^1(U^f;\Q)$ is semisimple. If $N$ is chosen such that $t^N=1$ on $H^1(U^f;\Q)$, let
\[ H^1(U^f;\Q)_1:= \ker \left( H^1(U^f;\Q) \xrightarrow{\cdot (t-1)} H^1(U^f;\Q) \right) \]
\[ H^1(U^f;\Q)_{\neq 1} := \ker \left( H^1(U^f;\Q) \xrightarrow{\cdot (t^{N-1} + \ldots + t+1)} H^1(U^f;\Q) \right). \]
Then we have an isomorphism of mixed Hodge structures \[ H^1(U^f;\Q)\cong H^1(U^f;\Q)_1\oplus H^1(U^f;\Q)_{\neq 1},\] 
and, moreover,  the following result holds (see \cite[Theorem 10.1.5]{EGHMW}):
\begin{thm}\label{thm:hyperplanes}
Let $\cA$ be an essential arrangement of $d$ hyperplanes in $\C^n$ defined by the zeros of a reduced polynomial $f$ of degree $d$, for $n\geq 2$. Then,
\begin{enumerate}
\item[(i)] $H^1(U^f;\Q)_1$ is a pure Hodge structure of type $(1,1)$, and has dimension $d-1$.
\item[(ii)] $H^1(U^f;\Q)_{\neq 1}$ is a pure Hodge structure of weight $1$.
\end{enumerate}
\end{thm}

To prove the above result, one first notices that, by a Lefschetz type argument, we can assume that $\cA$ is an essential line arrangent in $\C^2$. 
By the cohomological version of Theorem \ref{geoIntro}, the map \[ H^1(\pi)\colon H^1(U;\Q) \to H^1(U^f;\Q)\] is a mixed Hodge structure morphism, and Milnor's long exact sequence for $\pi\colon U^f \to U$ yields that  ${\rm Image} (H^1(\pi))=H^1(U^f;\Q)_1$. Part (i) of Theorem \ref{thm:hyperplanes} is then a consequence of the classical fact that  $H^1(U;\Q)$ is a pure Hodge structure of type $(1,1)$, see \cite{shapiro}. For part (ii), we use the cohomological version of Corollary~\ref{cor:fiber}, which yields a monomorphism of mixed Hodge structures \[H^1(U^f;\Q) \hookrightarrow H^1(F;\Q),\]
where $F$ is the generic fiber of $f\colon U \to \C^*$. 
The assertion follows by a careful analysis of the dimensions of the weight filtration on $H^1(F;\Q)$ (the only possible weights being $1$ and $2$). In fact, by \cite[Theorem 2.1]{DimcaTame} (and the discussion following it), the generic fiber of $f$ is connected, and one can show by direct computation that $\dim \Gr^W_2 H^1(F;\Q) = d-1$ (see \cite[Lemma 10.1.8]{EGHMW}).

\begin{remk} If $\cA$ is a \textit{central} hyperplane arrangement ($f$ is a homogeneous polynomial), then $f$ determines a global Milnor fibration with fiber $F$, so $H^{j}(U^f;\Q)\cong H^j(F;\Q)$ is an isomorphism of mixed Hodge structures for all $j$ (see Example \ref{ex4}).  Moreover, in this case the $t$-action is semisimple. Theorem \ref{thm:hyperplanes} provides a generalization (for $j=1$) of a similar result  for central arrangements (see, e.g., \cite[Theorem 7.7]{dimca2017hyperplanes} and the references therein). 
\end{remk}

\bibliographystyle{plain}

\begin{thebibliography}{10}


\bibitem{BBD} A.~A. Be\u{\i}linson, J.~Bernstein, and P.~Deligne.
\newblock Faisceaux pervers.
\newblock In {\em Analysis and topology on singular spaces, {I} ({L}uminy,  1981)}, volume 100 of {\em Ast\'{e}risque}, pages 5--171. Soc. Math. France,  Paris, 1982.


\bibitem{BudurLiuWang} N. Budur, Y. Liu, and B. Wang.
\newblock The monodromy theorem for compact {K}\"{a}hler manifolds and smooth quasi-projective varieties.
\newblock {\em Math. Ann.}, 371(3-4):1069--1086, 2018.





\bibitem{De} P.~Deligne.
\newblock Th\'{e}or\`eme de {L}efschetz et crit\`eres de  d\'{e}g\'{e}n\'{e}rescence de suites spectrales.
\newblock {\em Inst. Hautes \'{E}tudes Sci. Publ. Math.}, (35):259--278, 1968.

\bibitem{De2}
P.~Deligne.
\newblock Th\'{e}orie de {H}odge. {II}.
\newblock {\em Inst. Hautes \'{E}tudes Sci. Publ. Math.}, (40):5--57, 1971.



\bibitem{dimca1992hypersurfaces}
A.~Dimca.
\newblock {\em Singularities and topology of hypersurfaces}.
\newblock Universitext. Springer-Verlag, New York, 1992.


\bibitem{DimcaTame}
A.~Dimca.
\newblock Hyperplane arrangements, {$M$}-tame polynomials and twisted
  cohomology.
\newblock In {\em Commutative algebra, singularities and computer algebra
  ({S}inaia, 2002)}, volume 115 of {\em NATO Sci. Ser. II Math. Phys. Chem.},
  pages 113--126. Kluwer Acad. Publ., Dordrecht, 2003.


\bibitem{dimca2017hyperplanes}
A.~Dimca.
\newblock {\em Hyperplane arrangements. An introduction}.
\newblock Universitext. Springer, Cham, 2017.


\bibitem{DL}
A.~Dimca and A.~Libgober.
\newblock Regular functions transversal at infinity.
\newblock {\em Tohoku Math. J. (2)}, 58(4):549--564, 2006.

\bibitem{DN2}
A.~Dimca and A.~N\'{e}methi.
\newblock On the monodromy of complex polynomials.
\newblock {\em Duke Math. J.}, 108(2):199--209, 2001.


\bibitem{EGHMW} E.~Elduque, C.~Geske, M.~Herrad\'on Cueto, L.~ Maxim  and B.~Wang.
\newblock Mixed Hodge structures on Alexander modules.
\newblock arXiv:2002.01589.



\bibitem{thesiseva}
E.~Elduque.
\newblock Twisted Alexander modules of hyperplane arrangement complements, arXiv:1702.06267,
  2017.






\bibitem{hain1987rham}
R.~Hain.
\newblock The de {R}ham homotopy theory of complex algebraic varieties {I}.
\newblock {\em K-theory}, 1(3):271--324, 1987.








\bibitem{KK}
Vik.~S. Kulikov and V.~S. Kulikov.
\newblock On the monodromy and mixed {H}odge structure in the cohomology of an
  infinite cyclic covering of the complement to a plane curve.
\newblock {\em Izv. Ross. Akad. Nauk Ser. Mat.}, 59(2):143--162, 1995.

\bibitem{Lib94}
A.~Libgober.
\newblock Homotopy groups of the complements to singular hypersurfaces. {II}.
\newblock {\em Ann. of Math. (2)}, 139(1):117--144, 1994.

\bibitem{Lib96}
A.~Libgober.
\newblock Position of singularities of hypersurfaces and the topology of their
  complements.
\newblock J. Math. Sci.~82 (1): 3194--3210. 1996.

\bibitem{Liu}
Y.~Liu.
\newblock Nearby cycles and {A}lexander modules of hypersurface complements.
\newblock {\em Adv. Math.}, 291:330--361, 2016.




\bibitem{Max06}
L.~Maxim.
\newblock Intersection homology and {A}lexander modules of hypersurface
  complements.
\newblock {\em Comment. Math. Helv.}, 81(1):123--155, 2006.


\bibitem{MaxTommy}
L.~Maxim and K.~Wong.
\newblock Twisted {A}lexander invariants of complex hypersurface complements.
\newblock {\em Proc. Roy. Soc. Edinburgh Sect. A}, 148(5):1049--1073, 2018.

\bibitem{mil} J.~W. Milnor.
\newblock Singular points of complex hypersurfaces. 
\newblock Annals of Mathematics Studies, No. 61. Princeton University Press, Princeton, NJ; University of Tokyo Press, Tokyo (1968).



\bibitem{peters2008mixed}
C.~Peters and J.~Steenbrink.
\newblock {\em Mixed Hodge structures}, volume~52.
\newblock Springer Science \& Business Media, 2008.


\bibitem{shapiro}
B.~Z. Shapiro.
\newblock The mixed {H}odge structure of the complement to an arbitrary
  arrangement of affine complex hyperplanes is pure.
\newblock {\em Proc. Amer. Math. Soc.}, 117(4):931--933, 1993.







\end{thebibliography}

\end{document}